# Exploring Ring Structures: Multiset Dimension Analysis in Compressed Zero-Divisor Graphs

Nasir Ali[a,*], Hafiz Muhammad Afzal Siddiqui[a,] Muhammad Imran Qureshi[b]

[a]Department of Mathematics, COMSATS University Islamabad, Lahore Campus, Pakistan.

[b]Department of Mathematics, COMSATS University Islamabad, Vehari Campus, Pakistan.

**Email address:**

nasirzawar@gmail.com (Nasir Ali), hmasiddiqui@gmail.com (Hafiz Muhammad Afzal Siddiqui), imranqureshi18@gmail.com (Muhammad Imran Qureshi)

**Abstract:**

This paper explores the concept of multiset dimensions (*Mdim*) of compressed zero-divisor graphs (CZDG) associated with rings. The authors investigate the interplay between the ring-theoretic properties of a ring $R$ and the associated compressed zero-divisor graph. An undirected graph consisting of a vertex set $Z(R_E)\setminus\{[0]\} = R_E\setminus\{[0],[1]\}$, where $R_E = \{[x] : x \in R\}$ and $[x] = \{y \in R : ann(x) = ann(y)\}$ is called a compressed zero-divisor graph, denoted by $\Gamma_E(R)$. An edge is formed between two vertices $[x]$ and $[y]$ of $Z(R_E)$ if and only if $[x][y] = [xy] = [0]$, that is, iff $xy = 0$. For a ring $R$, graph $G$ is said to be realizable as $\Gamma_E(R)$ if $G$ is isomorphic to $\Gamma_E(R)$. We classify the rings based on *Mdim* of their associated CZDG and obtain the bounds for the *Mdim* of the compressed zero-divisor graphs. We also study the *Mdim* of realizable graphs of rings. Moreover, some examples are provided to support our results. Lately, we have discussed the interconnection between *Mdim*, girth, and diameter of CZDG.

## 1. Introduction

Graph theory has various applications in different fields. It is used to model interactions between individuals in social networks and to optimize routes in transportation systems. Graph analysis is beneficial for computer networks as it ensures data flow and connectivity. Electrical circuits can be better understood through graph representations, which aid in their design. In biology, graphs are used to depict protein interactions and genetic patterns. Epidemiology utilizes graphs to track the spread of diseases, while recommendation systems use them to suggest products online. Graphs are also used to model molecular structures in chemistry and to optimize search engines [13]. In strategic scenarios, game theory benefits from graph analysis [17].

Graph theory also intersects with algebra, leading to the development of algebraic graph theory. This area of study investigates the relationships between graphs and algebraic structures such as groups and matrices. Cayley graphs, for example, provide insights into the symmetries of groups. Spectral graph theory explores graph properties using eigenvalues and eigenvectors. Combinatorial optimization addresses problems such as maximal cliques and minimal spanning trees. Algebraic techniques, including graph theory, are helpful in designing error-correcting codes. Additionally, polynomials can be represented through graph-based interpolation, and readers may see relevant concepts [31, 32, 33].

Coding theory combines algebraic codes and graphs to identify and fix errors in data transmission. Representation theory investigates the connection between algebraic structures and graphs. Homological algebra examines the homology and cohomology of algebraic structures using graphs. Algebraic geometry benefits from graph representation, allowing for a visual understanding of algebraic varieties. The interplay between graph theory and algebra is commutative ringucial in various theoretical and practical situations.

Beck [7] proposed the connection between graph theory and algebra by introducing a Zero-divisor graph (ZD-graph) of a commutative ring $R$. The author's [7] primary focus was on the coloring of nodes in a graph, specifically on the ring elements that corresponded to these nodes. Note that a zero vertex is linked to all other vertices in this case. Let the set of zero divisors (ZD) is denoted by $Z(R)$ for a commutative ring, and the set of non-zero ZD of a commutative ring with $1 \neq 0$ is denoted by $Z^*(R) = Z(R)\backslash\{0\}$. In [3], Anderson and Livingston conducted a study on a ZD-graph in which each node represents a nonzero ZD. Let $x, y \in Z^*(R)$, then an undirected graph obtained by considering $x$ and $y$ as vertices forming an edge iff $xy = 0$ is called a ZD-graph of $R$, denoted $\Gamma(R)$. The study of Anderson and Livingston emphasizes the case of finite rings, as finite graphs can be obtained when $R$ is finite. Their task was to determine whether a graph is complete for a given ring or a star for a given ring. This ZD-graph definition differs slightly from Beck's ZD-graph definition for $R$. Remember, zero is not considered as a vertex of the ZD-graph in this case. The study of ZD-graph has been extinct in recent years, and the idea has been explored, which leads us to the new form of ZD-graph that includes ideal-based ZD-graph and module-based ZD-graph [15, 20, 23]. Redmond [22] expanded the ZD-graph idea from unital commutative rings to noncommutative rings. Different methods were presented by him to characterize the ZD-graph related to a noncommutative ring, encompassing both undirected & directed graphs. Redmond extended this work using a ZD-graph for a commutative ring and transformed it into an ideal-based ZD-graph. The aim was to generalize the method by substituting elements with zero products with elements whose product belongs to a particular ideal $I$ of ring $R$.

Mulay's [16] work inspired us to study the ZD-graph obtained by considering equivalence classes of ZD of a ring $R$. This type of ZD-graph is called CZDG, denoted by $\Gamma_E(R)$ [4]. A CZDG is an undirected graph obtained by considering $Z(R_E)\backslash\{[0]\} = R_E\backslash\{[0],[1]\}$ as vertex set, and can be constructed by taking the equivalence classes $[x] = \{y \in R : ann(x) = ann(y)\}$, for every $x \in R\backslash([0]\cup[1])$ as vertices and an edge is formed between two distinct classes $[x]$ and $[y]$ iff $[x][y] = 0$, i.e., iff $xy = 0$. It is important to note that if two vertices say $x$ and $y$, are adjacent in $\Gamma(R)$, then in CZDG, $[x]$ and $[y]$ are adjacent iff $[x] \neq [y]$. Clearly, $[1] = R\backslash Z(R)$ and $[0] = \{0\}$, also for each $x \in R \setminus ([0] \cup [1])$, $[x] \subseteq Z(R)\backslash\{0\}$. Readers may study [5] for some interesting results on CZDG.

We consider an example to understand the concept of ZD-graph and CZDG. Let $R = \mathbb{Z}_{16}$, then the vertex set of $\Gamma(R) = \{2, 4, 6, 8, 10, 12, 14\}$, see Figure 1 (i) shows its ZD-graph. Now, we see $ann(2) = \{8\}$, $ann(4) = \{4, 8, 12\}$, $ann(6) = \{8\}$, $ann(8) = \{2, 4, 6, 8, 10, 12, 14\}$, $ann(10) = \{8\}$, $ann(12) = \{4,8,12\}$ and $ann(14) = \{6,8\}$. Hence, the vertex set for $\Gamma_E(R) = \{[2], [4], [8], [14]\}$. See Figure 1 (ii) for its CZDG.

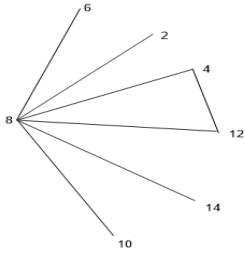 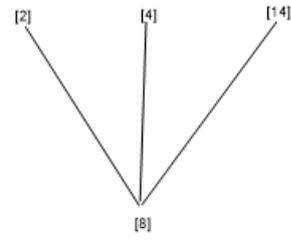

Figure 1(i)    Figure 1(ii)

The Annihilator ideals in the ring correspond to the vertices of the $\Gamma_E(R)$. Also, remember that $diam(\Gamma_E(R)) \leq 2$, and CZDG is always connected. Also, $diam(\Gamma_E(R)) \leq diam(\Gamma(R))$. For the CZDG, $gr(\Gamma_E(R)) \leq 3$ whenever CZDG of $R$ contains a cycle [5]. It can be seen in [4] that $\Gamma_E(R)=\Gamma_E(S)$ if $S$ is a Noetherian or $S$ is a finite commutative ring.

Readers may see [1, 2, 3] to read many advantages of studying CZDG over the earlier studied ZD-graph. For example, in any ring $R$ having at least 2 vertices, there exists no finite regular CZDG [[27], Proposition 1.10]. Further, Spiroff et al. [27] showed that the CZDG of local ring $R$ is isomorphic to a star graph with a minimum of 4 vertices. (If a ring $R$ has a unique maximal ideal, then it is called a local ring).

Another commutative ringucial aspect to consider is the connection between studying equivalence class graphs. The associated primes of $R$ usually considered distinct vertices in CZDG. In this paper, all graphs are simple graphs; a commutative ring with unity is denoted by $R$, and the units set is considered as $U(R)$. $\mathbb{Z}_n$ denotes ring of integers modulo $n$ and $\mathbb{F}_q$ denotes finite field on $q$ elements. Readers are encouraged to study [11, 17, 31] for basic definitions of graph theory and [6, 13] to study basic definitions of ring theory.

Formally, the graph is an ordered pair $G = (V, E)$; here, $V$ and $E$ denotes vertices and edge set, respectively. A graph's order and size are defined as the cardinality of nodes and edges set, respectively. The open neighborhood of a node $v$ is written as $N(v)$, and defined as $\{v \in V(G) : vu \in E(G)\}$, while the closed neighborhood of a node $u$ is written as $N[u]$, and defined as $\{u\} \cup N(u)$. The distance between two nodes $u'$ and $v'$ is denoted by $d(u', v')$ and defined as the length of the shortest path between them, while $d(w, e') = min\{d(w, u'), d(w, v')\}$ defines the distance between a node $w$ and the edge $e' = u'v'$.

A graph $G$ is a regular graph if for every $r \in V, deg(r) = c$ for a fixed $c \in Z^+$. A complete graph is a graph in which there is a connection between every pair of vertices, and we denote it by $k_m$, where $m$ stands for the number of vertices. A graph is classified as a complete bipartite graph if we can divide its vertices into 2 distinct sets, $X$ and $Y$, where every node in $X$ is connected to every node in $Y$, and it is usually denoted by $k_{m,n}$, where $|X| = m$ and $|Y| = n$. In a connected graph $G$, when the removal of a vertex results in the formation of two or more components, then it is called a cut vertex.

In [30], the authors have studied the edge metric dimension (EMD) of various graphs. Moreover, the relationship between the metric dimension (MD) and EMD allows for the identification of graphs where these two dimensions are equal, as well as for some other graphs $G$ for which $dim(G) < dim_E(G)$ or $dim_E(G) < dim(G)$. Basically, Kelenc et al. explored the comparison of values $dim(G)$ and $dim_E(G)$. Recently, a study on metric parameters for ZD-graphs has been done. Redmond, in 2003 [23], studied the ideal-based ZD-graph of commutative rings was studied by him.

Simanjuntak et al. [28] introduced a new variant of metric dimension known as multiset dimension (*Mdim*), where distances between $v$ and all vertices in resolving set (RS) $W$ were calculated, including their multiplicities. The *Mdim* is defined as the minimum cardinality of the RS. Assume $v$ is a node of $G$ and $B \subseteq V(G)$, then representation for multiset of $v$ w.r.t $B$ defined as the distances between $v$ and nodes in $B$. This representation is denoted by $r_m(B)$. For all pairs of distinct nodes $u$ and $v$, $B$ is called $m-$resolving set of $G$, if $r_m(B) \neq r_m(B)$. The cardinality of $m-$resolving set is called the multiset basis of the graph $G$, and the minimum cardinality of the multiset basis is called the multiset dimension of $G$; we denote it by $Mdim(G)$, if $G$ does not contain a $m-$resolving set, we write $Mdim(G) = \infty$. The key point of this article is that the apparent expansions are an oversimplification of the task of identifying graph vertices using the multiset representation.

Some of the contributions of this article include generalizing the *Mdim* of CZDG that accommodates different characterizations of the rings based on vertices in CZDG. Authors have also proved that rings can be characterized based on their multiset dimensions. Moreover, authors have characterized several rings based on their CZDG for which multiset dimension can be bounded by diameter.

The novelty of finding the multiset dimensions of graphs lies in the fact that it is a relatively new concept and has not been extensively studied. It provides a more complete understanding of the graph's structure and algebraic properties, which can be useful in many applications such as network design, social networking, and communication systems.

Simanjuntak et al. [28] found some sharp bounds for *Mdim* of arbitrary graphs in terms of its MD, order, or diameter. For a graph to have finite *Mdim*, Siamanjuntak provided some necessary conditions, with an example of an infinite family of graphs where those necessary conditions are also sufficient. It was also shown that the *Mdim* of any graph other than a path is at least 3, and two families of graphs having the *Mdim* 3 were proved. Here, we consider some results from [28] as follows:

**Theorem 1.1 [28]:** For any integer $m \geq 3, n \geq 6$, $Mdim(P_m) = 1$, $Mdim(k_m) = \infty$. and $Mdim(C_n) = 3$. Moreover, $Mdim(P_m) = 1$ iff $G \cong P_m$.

Moreover, for a complete bipartite graph $k_{m,n}$, we get different *Mdim* for different choices of values of $m$ & $n$.

**Theorem 1.2 [28]:** For any complete bipartite graph $k_{m,n}$, the $Mdim(k_{m,n})$ is given below.

| | |
|---|---|
| $Mdim(k_{m,n}) = 1$ | $for\ m = 1\ and\ n = 1,2$ |
| | $for\ m = 2\ and\ n = 1$ |
| $Mdim(k_{m,n}) = \infty$ | $for\ m = 1\ and\ n \geq 3$ |
| | $for\ m = 2\ and\ n \geq 2$ |

Moreover, *Mdim* for a single vertex graph $G$ is supposed to be zero and undefined for an empty graph.

## 2. Multiset dimension of some $\Gamma_E(R)$

**Proposition 2.1:** Let $R$ be a finite commutative ring. Then $Mdim(\Gamma_E(R)) = 0$ iff $\Gamma(R) \cong k_n$.

**Proof.** Let $\Gamma(R) \cong k_n$, then one of the two cases holds, either $R$ is isomorphic to $\mathbb{Z}_2 \times \mathbb{Z}_2$ or for all $x, y \in Z^*(R)$ product of $x$ and $y$ is zero. i.e., $xy = 0$. Let $v_1, v_2, \ldots v_n$ be the ZD of the ring $R$ that corresponds to the vertices of the ZD-graph, then $[v_1] = [v_2], \cdots = [v_n]$ suggests that all vertices of the ZD-graph would collapse to a graph with a single vertex in CZDG, and as we know, $Mdim$ is 0 for a single vertex graph.

For the converse part, assume that $\Gamma(R) \not\cong k_n$. This suggests the presence of at least one vertex in $\Gamma(R)$ that is not connected to every other vertex. Consequently, the size of the edge set CZDG is at least 2, i.e., $|\Gamma_E(R)| \geq 2$, leading to a multiset dimension $Mdim(\Gamma_E(R)) \geq 1$.

The converse can also be shown by assuming that $Mdim(\Gamma_E(R)) = 0$. This implies that $\Gamma_E(R) = \{[a]\}$ for some non-zero element say, $a \in Z^*(R)$, indicating that $\Gamma_E(R)$ is a graph with single vertex. This implies that $\Gamma(R)$ is either isomorphic to $k_n, \forall, n \geq 1$ or a graph with a single vertex. □

**Proposition 2.2:** Let $R$ be a finite commutative ring. Then $Mdim(\Gamma_E(R)) = 1$, iff $\Gamma(R) \cong k_{m,n}$, where $m$ or $n \geq 2$.

**Proof.** Suppose that the ZD-graph $\Gamma(R)$ is isomorphic to $k_{m,n}$ having two distance similar classes $V_1$ and $V_2$. Specifically, assume that $V_1 = \{u_1, u_2, \cdots, u_m\}$ and $V_2 = \{v_1, v_2, \cdots, v_m\}$ such that $u_i v_j = 0, \forall i \neq j$. Cleary, an independent set is formed by each of $V_1$ and $V_2$. Furthermore, observing that $[u_1] = [u_2] = \cdots = [u_m]$ and also, $[v_1] = [v_2] = \cdots = [v_n]$, we deduce that both $V_1$ and $V_2$ each represents a single vertex CZDG. Given the connected nature of the graph, we can conclude that $\Gamma_E(R) \cong k_{1,1}$, which can be visualized as a path consisting of two vertices. Therefore, by Theorem 1.2, it follows that the multiset dimension of CZDG is 1 that is $Mdim(\Gamma_E(R)) = 1$. □

**Remark 2.1:** It's important to note that the converse of Proposition 2.2 may not hold; this is exemplified by the graph illustrated in Figure 2.1. However, in cases where $R \cong \mathbb{Z}_2 \times \mathbb{Z}_2$, we find that $\Gamma_E(R) \cong k_{1,1}$, with $Mdim(\Gamma_E(R)) = 1$ such that $\Gamma(R) \cong \Gamma_E(R)$. As demonstrated in ([27], Proposition 1.5), a significant difference between ZD-graph and CZDG is that the latter cannot be isomorphic to a graph having a minimum of 3 nodes and is also complete. However, see ([27], Proposition 1.7), if $\Gamma_E(R)$ is isomorphic to a complete r-partite graph, then $r$ must equal 2, resulting in $\Gamma_E(R) \cong k_{n,1}$, for some $n \geq 1$. A ring R is classified as a Boolean ring if $a^2 = a$ holds for every element $a \in R$. Importantly, a Boolean ring $R$ is both commutative and has a characteristic of 2 ($char(R) = 2$). More comprehensively, a commutative ring qualifies as a von Neumann regular ring. When, for any element $a \in R$, there exists an element $b$ within $R$ such that $a = a^2 b$. This condition is equivalent to $R$ being a zero-dimensional reduced ring, as elucidated in ([13], Theorem 3.1). Clearly, a Boolean ring can be classify as a von Neumann regular ring, but the converse may not always be true. For instance, consider a family $\{F_i\}_{i \in I}$ of fields, where the product $\prod_{i \in I} F_i$ is von Neumann regular ring. Nevertheless, it is Boolean if and only if $F_i \cong \mathbb{Z}_2$ holds for all $i \in I$. If $R$ is reduced ring, then for $r, s \in \Gamma(R)$, the conditions $N(r) = N(s)$ and $[r] = [s]$ are equivalent ([13], Lemma 3.1). Furthermore, if $R$ is a von Neumann regular ring, then these conditions are equivalent to $rR = sR$.

Moreover, if $R$ is a von Neumann regular ring and $B(R)$ signifies the set of idempotent elements within $R$, the mapping defined by $e \mapsto [e]$ forms an isomorphism from the subgraph of $\Gamma(R)$ induced by $B(R)\setminus\{0,1\}$ onto $\Gamma_E(R)$ ([13], Proposition 4.5). Particularly, if $R$ is a Boolean ring (i.e., $R = B(R)$), then $\Gamma_E(R)$ is isomorphic to $\Gamma(R)$. This discourse leads us to the subsequent characterization.

**Corollary 2.1:** Consider $R$ and $S$ bring reduced commutative ring. If $\Gamma(R)$ is isomorphic to $\Gamma(S)$ then $Mdim(\Gamma_E(R)) = Mdim(\Gamma_E(S))$.

## 3. Analyzing bounds for the multiset dimension of $\Gamma_E(R)$

Here, we explore the importance of calculating *Mdim* in studying CZDG. Additionally, we determine the *Mdim* of certain specialized ring types that correspond to $\Gamma_E(R)$. Notably, a recent contribution by Pirzada et al. [18] presented a work on the characterization of $\Gamma(R)$ in cases where the md is finite and cases where it remains undefined ([18], Theorem 3.1).

**Theorem 3.1:** Let R be a finite commutative ring with unity. Then $Mdim(\Gamma_E(R))$ is undefined iff $R$ is an integral domain.

**Proof.** As we know that $\Gamma_E(R)$ is not defined if $R$ is an integral domain. Which follows that $Mdim(\Gamma_E(R))$ is undefined and vice versa. □

Let us consider the following lemma,

**Lemma 3.1:** Consider a ring $R$, which is a finite local ring, $|R| = p^n$, with $n$ being a positive integer and some prime $p$.

The preceding result will be applied to determine the *Mdim* of rings $R$, which are finite local rings. □

**Proposition 3.1:** Let $|R| = p^2$ for a local ring $R$ and $p$ being 2, 3 or 5, then $Mdim(\Gamma_E(R))$ can be either undefined or 0.

**Proof.** Let us assume all local rings have an order. $p^2$, where $p$ is a prime. The following rings $\mathbb{F}_{p^2}$, $\mathbb{Z}_{p^2}$ and $\frac{\mathbb{F}_p[x]}{(x^2)}$ are local rings of order $p^2$ ([10], p. 687).

Case I: If $R \cong \mathbb{F}_{p^2}$, i.e., $R$ is a field having an order $p^2$. In such cases, the graph $\Gamma_E(R)$ becomes an empty graph, resulting in an undefined $Mdim(\Gamma_E(R))$.

Case II: If $R \cong \mathbb{Z}_{p^2}$, or $R \cong \frac{\mathbb{F}_p[x]}{(x^2)}$, i.e., $|R| = p^2$ and $R$ is not a field. For $p = 2, 3,$ or 5, the graph $\Gamma_E(R)$ consists of a single vertex, leading to $Mdim(\Gamma_E(R)) = 0$. This concludes our result. □

**Example 3.1:** Consider the ring $R \cong \mathbb{Z}_8$, then it is easy to see that $\Gamma(R) \cong P_3$ $and$ $\Gamma_E(R) \cong P_2$. Hence by Theorem 1.1, $\text{Mdim}(\Gamma(R)) = 1 = Mdim(\Gamma_E(R))$.

**Proposition 3.2:** Let $R$ be a local ring having order (i.e., $R$ is not a field)

(i) $p^3$ with $p = 2$ $or$ $3$, then $Mdim(\Gamma_E(R)) = 0$, and $Mdim(\Gamma_E(R)) = 1$ only if $R \cong \mathbb{Z}_8$, $\mathbb{Z}_{27}, \mathbb{Z}_2[x]/(x^3)$, $\mathbb{Z}_3[x]/(x^3)$, $\mathbb{Z}_4[x]/(2x, x^2 - 2)$, $\mathbb{Z}_9[x]/(3x, x^2 - 3)$, $\mathbb{Z}_9[x]/(3x, x^2 - 6)$.

(ii) $p^4$ with $p = 2$, then $Mdim(\Gamma_E(R))$ is equal to 0, 1, or $\infty$.

**Proof.** (a) Following rings $\mathbb{F}_{p^3}$, $\frac{\mathbb{F}_p[x]}{(x^3)}, \frac{\mathbb{F}_p[x,y]}{(x,y)^2}, \frac{\mathbb{Z}_{p^2}[x]}{(px, \ x^2)}, \frac{\mathbb{Z}_{p^2}[x]}{(px, \ x^2-p)}$ are all local rings having an order $p^3$.

Case (i). When $p = 2$, the local rings $\frac{\mathbb{Z}_2[x,y]}{(x,y)^2}$ and $\frac{\mathbb{Z}_4[x]}{(2x, \ x^2)}$ have the same equivalence classes of zero divisors. Let us assume that $[m] = \{x, y, x + y\}$ is equivalence class for any ZD $m$ of ring $\frac{\mathbb{Z}_2[x,y]}{(x,y)^2}$ and $[n] = \{2, x, x + 2\}$ is equivalence class for any ZD $m$ ring $\frac{\mathbb{Z}_4[x]}{(2x, \ x^2)}$, i.e., $\Gamma_E(R)$ is a single vertex graph. Hence, $Mdim(\Gamma_E(R)) = 0$. But $\Gamma_E(R) \cong P_2$ for $\mathbb{Z}_8$, $\mathbb{Z}_2[x]/(x^3)$, $\mathbb{Z}_4[x]/(2x, x^2 - 2)$. Hence by Theorem 1.1 $Mdim(\Gamma_E(R)) = 1$.

Case (ii). If $p = 3$ for the above given rings, we found that the CZDG structures of the rings $\mathbb{Z}_{27}$, $\mathbb{Z}_3[x]/(x^3)$, $\mathbb{Z}_9[x]/(3x, x^2 - 3)$, $\mathbb{Z}_9[x]/(3x, x^2 - 6)$ is the same and isomorphic to $P_2$. Then by Theorem 1.1 $Mdim(\Gamma_E(R)) = 1$.

Also, the rings $\frac{\mathbb{Z}_{3^2}[x]}{(3x,x^2)}$ and $\frac{\mathbb{Z}_3[x,y]}{(x,y)^2}$ have same equivalence classes of ZD given by $[m] = \{3, 6, x, \ 2x, \ x + 3, x + 6, 2x + 3, 2x + 6\}$ for any ZD $m$ of ring $\frac{\mathbb{Z}_{3^2}[x]}{(3x,x^2)}$ and $[n] = \{x, \ 2x, \ y, \ 2y, \ x + y, \ 2x + y, \ x + 2y, \ 2x + 2y\}$ for any ZD $n$ of ring $\frac{\mathbb{Z}_3[x,y]}{(x,y)^2}$, i.e., $\Gamma_E(R)$ is a graph with a single vertex. Hence, $Mdim(\Gamma_E(R)) = 0$.

(ii) Now we focus our attention to local rings having order $p^4$, by taking $p = 2$. It is found in [10] that there are 21 non-isomorphic commutative local rings with an identity of order 16. The following are rings with $Mdim(\Gamma_E(R)) = 0$ are $\mathbb{F}_4[x]/(x^2)$, $\mathbb{Z}_2[x, y, z]/(x, y, z)^2$, and $\mathbb{Z}_4[x]/(x^2 + x + 1)$. The following are the rings having $Mdim(\Gamma_E(R)) = 1$, $\mathbb{Z}_2[x]/(x^4)$, $\mathbb{Z}_2[x, y]/(x^3, xy, x^2)$, $\mathbb{Z}_4[x]/(2x, x^3 - 2)$, $\mathbb{Z}_4[x]/(x^2 - 2), \mathbb{Z}_8[x]/(2x, \ x^2), \mathbb{Z}_{16}, \mathbb{Z}_4[x]/(x^2 - 2x - 2), \mathbb{Z}_8[x]/(2x, \ x^2 - 2)$ and $\mathbb{Z}_4[x]/(x^2 - 2x)$. Furthermore, the following rings $\mathbb{Z}_4[x]/(x^2)$, $\mathbb{Z}_2[x, y]/(x^2, y^2)$ and $\mathbb{Z}_2[x, y]/(x^2 - y^2, xy)$ have $Mdim(\Gamma_E(R)) = \infty$. $\square$

Now, we will find $Mdim(\Gamma_E(\mathbb{Z}_n))$.

**Proposition 3.3:** Consider a prime number $p$.
(a) When $n = 2p$ & $p > 2$, then $Mdim(\Gamma_E(\mathbb{Z}_n)) = 1$.

(b) When $n = p^2$, then $Mdim(\Gamma_E(\mathbb{Z}_n)) = 0$.

**Proof.** (a) When $p = 2$, CZDG of $\mathbb{Z}_4$ is a single vertex graph, hence $Mdim(\Gamma_E(\mathbb{Z}_4)) = 0$.

Consider $p > 2$, then $\{2, 2.2, 2.3, \ldots, 2.(p-1), p\}$ is ZD set of $\mathbb{Z}_n$. Given that the $Char\ (\mathbb{Z}_n) = 2p$, we can deduce that $p$ is adjacent to every other vertex. Consequently, the equivalence classes of these ZD are as follows:
$$[2] = [2.2] = \ldots = [2.(p-1)] = \{p\}, \qquad [p] = \{2, 2.2, 2.3, \ldots, 2.(P-1)\}.$$
Thus, $\{[p], [2x]\}$ is vertex set of $\Gamma_E(\mathbb{Z}_n)$, for some positive integer $x = 1, 2, 3, \ldots, p-1$. So $\Gamma_E(\mathbb{Z}_n) \cong P_2$ then, by Theorem 1.1 $Mdim(\Gamma_E(R)) = 1$.

(b) ZD set of $\mathbb{Z}_n$ is $\{p, p.2, p.3, \ldots, p.(p-1)\}$, when $n = p^2$ and $p > 2$. Since $Char\ (\mathbb{Z}_n) = p^2$, so equivalence classes of all are the same, i.e., $\{p, p.2, p.3, \ldots, p.(p-1)\}$. Hence CZDG of $\mathbb{Z}_n$ is a single vertex graph, so $Mdim(\Gamma_E(R)) = 0$.

Consider a ring $R$, a graph $G$ is said to be realizable as $\Gamma_E(R)$ if $G \cong \Gamma_E(R)$. However, numerous results suggest that the majority of graphs cannot be realized as $\Gamma_E(R)$, For instance, $\Gamma_E(R)$ cannot represent a complete graph with three or more vertices or a cycle graph. □

**Proposition 3.4:** Consider a realizable graph $\Gamma_E(R)$ containing 3 vertices, then $Mdim(\Gamma_E(R)) = 1$.

**Proof.** In a study by Spiroff et al. [27], it was established that one realizable graph denoted by $\Gamma_E(R)$ encompassing precisely 3 vertices as a graph representing equivalence classes of zero-divisors for a given ring $R$ is isomorphic to the path graph $P_3$, see Figure 2. understandably, its $Mdim$ is 1. □

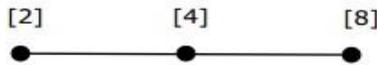

Figure 2: $\mathbb{Z}_{16}$

**Proposition 3.5:** Consider a realizable graph $\Gamma_E(R)$ containing 4 vertices, then $Mdim(\Gamma_E(R))$ is either $\infty$ or 1.

**Poof.** The possible realizable graph $\Gamma_E(R)$ with four vertices are illustrated in Figure 3. It is readily apparent that their $Mdim$ can be either $\infty$ or 1. □

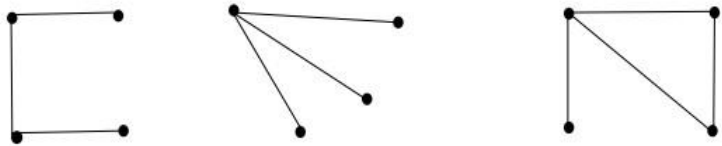

Figure 3: $(\mathbb{Z}_4 \times \mathbb{F}_4)$, $\mathbb{Z}_4[x]/(x^2)$, $\mathbb{Z}[x,y]/(x^3, xy)$

**Proposition 3.6:** Consider a realizable graph $\Gamma_E(R)$ with 5 vertices, then $Mdim(\Gamma_E(R)) = \infty$.

**Poof.** The possible realizable graph $\Gamma_E(R)$ with 5 vertices are illustrated in Figure 4; It is readily apparent that their *Mdim* is ∞. □

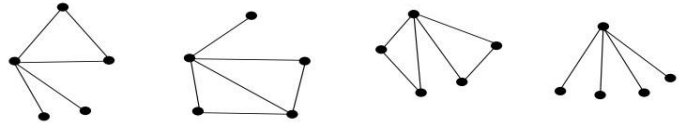

Figure 4: $\mathbb{Z}_9[x]/(x^2), \mathbb{Z}_{64}, \mathbb{Z}_3[x,y]/(xy, x^3, y^3, x^2 - y^2), \mathbb{Z}_8[x,y]/(x^2, y^2, 4x, \ 4y, \ 2xy)$

## 4. Relationship between multiset dimension, diameter, and girth of $\Gamma_E(R)$

Within this section, we delve into the correlation among multiset dimensions, diameter, and girth of $\Gamma_E(R)$. Since $gr(\Gamma_E(R)) \in \{3, \infty\}$, but $gr(\Gamma_E(R)) = 3$ iff $gr(\Gamma(R)) = 3$, where $R$ is a reduced commutative ring. Moreover, $gr(\Gamma_E(R)) = \infty \Leftrightarrow gr(\Gamma(R)) \in \{4, \infty\}$. Nevertheless, it is possible that $gr(\Gamma(R)) = 3$ and we could have either $gr(\Gamma_E(R)) = 3$ or $\infty$. In the subsequent result, we establish the $Mdim(\Gamma_E(R))$ of a ring $R$ in terms of $gr(\Gamma_E(R))$.

**Theorem 4.1:** Consider a finite commutative ring $R$ such that $gr(\Gamma_E(R)) = \infty$.
(i)   Then $Mdim(\Gamma_E(R)) = 1$, if $R$ is a reduced ring,
(ii)  Then $Mdim(\Gamma_E(R)) = 0$, if $R \cong \mathbb{Z}_4, \mathbb{Z}_9, \mathbb{Z}_2[x]/(x^2)$.

**Proof.** (a) Let us consider that $R$ is a reduced ring and not isomorphic to $\mathbb{Z}_2 \times \mathbb{Z}_2$, then it is certain that $R$ cannot be isomorphic to $\mathbb{Z}_2 \times A$, where $A$ is some finite field. Consequently, $R$ possesses two equivalence classes of ZD [(1,0)] & [(0,1)], which are adjacent to each other (by Theorem 2 [29]). As a result, $Mdim(\Gamma_E(R)) = 1$. However, if $R$ is isomorphic to $\mathbb{Z}_2 \times \mathbb{Z}_2$, suggesting that $R$ is a Boolean ring. It follows that $\Gamma(R)$ is isomorphic to $\Gamma_E(R)$. Consequently, the result is completed by using (Lemma 3, [29]).

(b) For the rings listed as follows: $R \cong \mathbb{Z}_4, \mathbb{Z}_9, \mathbb{Z}_2[x]/(x^2)$, their CZDG, $\Gamma_E(R)$ represents a single vertex graph and hence $Mdim(\Gamma_E(R)) = 0$. □

**Corollary 4.1:** Let $R$ be a finite commutative ring having unity 1 and $\mathbb{F}$ be a finite field. Also, let $R$ be a local ring isomorphic to any of the following listed rings, $\mathbb{Z}_8, \mathbb{Z}_{27}, \mathbb{Z}_2[x]/(x^3), \mathbb{Z}_4[x]/(2x, x^2 - 2), \mathbb{Z}_2[x,y]/(x^3, \ xy, \ y^2), \mathbb{Z}_8[x]/(2x, x^2), \mathbb{Z}_4[x]/(x^3, 2x^2, 2x), \mathbb{Z}_9[x]/(3x, x^2 - 6), \mathbb{Z}_9[x]/(3x, x^2 - 3), \mathbb{Z}_3[x]/(x^3)$. If $gr(\Gamma_E(R)) = \infty$, then for reduced rings $R \times \mathbb{F}$, $\Gamma_E(R \times \mathbb{F}) \cong \Gamma_E(R)$ with *Mdim* 1.

**Proof.** Let $gr(\Gamma_E(R)) = \infty$, then for all reduced rings $R \times \mathbb{F}$, we have $\Gamma_E(R \times \mathbb{F}) \cong k_{1,1}$ (Lemma 3, [29]). Furthermore, the above local rings list has the same CZDG isomorphic to $k_{1,1}$. Hence $Mdim(\Gamma_E(R)) = 1$ for the above local rings. □

We now focus on examining the connection between multiset dimension and diameter of CZDG. Considering, $diam\ (\Gamma_E(R)) \leq 3$, when there is a cycle in CZDG, the following outcomes emerge.

**Theorem 4.2:** Let $R$ be a commutative ring and $\mathbb{F}_1$ and $\mathbb{F}_2$ are fields. Then,

(a) $Mdim(\Gamma_E(R)) = diam(\Gamma_E(R)) = 1$ if $R \cong \mathbb{F}_1 \times \mathbb{F}_2$.

(b) $Mdim(\Gamma_E(R)) = 0 \Leftrightarrow diam(\Gamma_E(R)) = 0$.

(c) $Mdim(\Gamma_E(R)) = 0$ if $Z(R)^2 = 0$ and $|Z(R)| \geq 2$.

(d) $Mdim(\Gamma_E(R)) = 0 \Leftrightarrow diam(\Gamma(R)) = 0$ or $1, R \not\cong \mathbb{Z}_2 \times \mathbb{Z}_2$.

**Proof.**

(a) Let $R \cong \mathbb{F}_1 \times \mathbb{F}_2$, then $|\Gamma_E(R)| = 2$ (by Lemma 3, [29]). Since only equivalence classes of ZD are $[(0,1)]$ and $[(1,0)]$. So, CZDG is isomorphic to $k_{1,1}$. Hence, $Mdim(\Gamma_E(R)) = 1 = diam(\Gamma_E(R))$.

(b) $Mdim(\Gamma_E(R)) = 0$ iff $\Gamma_E(R)$ is a single vertex graph, and it is possible iff $diam(\Gamma_E(R)) = 0$.

(c) Let $|Z(R)| \geq 2$ and $Z(R)^2 = 0$. Hence $ann(a) = ann(b)$, for each $a, b \in Z(R)^*$, This suggests that $diam(\Gamma_E(R)) = 0$. Therefore $Mdim(\Gamma_E(R)) = 0$.

(d) Assume that $diam(\Gamma(R)) = 0$ or $1$, subsequently $\Gamma(R) \cong k_n$, thus $Mdim(\Gamma_E(R)) = 0$ except if $R \not\cong \mathbb{Z}_2 \times \mathbb{Z}_2$. Conversely, assume that $Mdim(\Gamma_E(R)) = 0$. Then $\Gamma(R) \cong k_n$, hence $diam(\Gamma(R)) = 0$ or $1$. □

## 5. Conclusion

The rings have been characterized by studying the multiset dimensions of their associated CZDG. We found some bounds for the multiset dimensions of CZDG. We also studied the *Mdim* of graph $G$, which is realizable as $\Gamma_E(R)$. The study concluded by discussing the connection between multiset dimension, diameter, and girth of $\Gamma_E(R)$. In the future, *Mdim* can be studied to characterize the rings based on their other types of ZD-graphs.

**Declaration:**

- Availability of data and materials: The data is provided on request to the authors.
- Conflicts of interest: The authors declare that they have no conflicts of interest, and all agree to publish this paper under academic ethics.
- Competing Interests and Funding: The authors did not receive support from any organization for the submitted work.
- Author's contribution: All the authors equally contributed to this work.
- Fundings: This received no specific grant from any funding agency in the public, commercial, or not-for-profit sectors.

**References:**

[1]  S. Akbari and A. Mohammadian, On zero divisor graphs of a ring, *J. Algebra*, 274 (2004), 847–855.

[2]  D. D. Anderson and M. Naseer, Beck's coloring of a commutative ring, *J. Algebra*, 159 (1993), 500–514.

[3]  D. F. Anderson and P. S. Livingston, The zero divisor graph of a commutative ring, *J. Algebra*, 217 (1999), 434–447.

[4]  D. F. Anderson and J. D. LaGrange, Commutative Boolean Monoids, reduced rings and the compressed zero-divisor graphs, *J. Pure and Appl. Algebra*, 216 (2012), 1626–1636.

[5]  D. F. Anderson and J.D. LaGrange, Some remarks on the compressed zero-divisor graphs, *J. Algebra*, 447 (2016), 297–321.

[6]  M. F. Atiyah and I. G. MacDonald, *Introduction to Commutative Algebra*, Addison-Wesley, Reading, MA (1969).

[7]  I. Beck, Coloring of commutative rings, *J. Algebra*, 26 (1988), 208–226.

[8]  P. J. Cameron and J. H. Van Lint, *Designs, Graphs, Codes, and their Links*, in: London Mathematical Society Student Texts, 24 (5), Cambridge University Press, Cambridge, 1991.

[9]  G. Chartrand, L. Eroh, M. A. Johnson and O. R. Oellermann, Resolvability in graphs and metric dimension of a graph, *Disc. Appl. Math.*, 105 (1-3) (2000), 99–113.

[10]  B. Corbas and G. D. Williams, Rings of order $p^5$. II. Local rings, *J. Algebra*, 231 (2) (2000), 691–704.

[11]  R. Diestel, Graph Theory, 4th ed. Vol. 173 of Graduate texts in mathematics, Heidelberg: Springer, 2010.

[12]  F. Harary and R. A. Melter, On the metric dimension of a graph, *Ars Combin.*, 2 (1976), 191–195.

[13]  J. A. Huckaba, *Commutative Rings with Zero Divisors*, Marcel-Dekker, New York, Basel, 1988.

[14]  S. Khuller, B. Raghavachari and A. Rosenfeld, *Localization in graphs*, Technical Report CS-TR-3326, University of Maryland at College Park, 1994.

[15]  H. R. Maimani, M. R. Pournaki and S. Yassemi, Zero divisor graph with respect to an ideal, *Comm. Algebra*, 34 (3) (2006), 923–929.

[16]  S. B. Mulay, Cycles and symmetries of zero-divisors, *Comm. Algebra*, 30 (7) (2002), 3533–3558.

[17]  S. Pirzada, An Introduction to Graph Theory, Universities Press, Orient Blackswan, Hyderabad, 2012.

[18]  S. Pirzada, Rameez Raja and S. P. Redmond, Locating sets and numbers of graphs associated to commutative rings, *J. Algebra Appl.*, 13 (7) (2014), 1450047.

[19]  S. Pirzada and Rameez Raja, On the metric dimension of a zero divisor graph, *Comm. Algebra*, 45 (4) (2017), 1399–1408.

[20]  S. Pirzada and Rameez Raja, On graphs associated with modules over commutative rings, *J. Korean Math. Soc.*, 53 (5) (2016), 1167–1182.

[21]  Rameez Raja, S. Pirzada and Shane Redmond, On locating numbers and codes of zero divisor graphs associated with commutative rings, *J. Algebra Appl.*, 15 (1) (2016), 1650014.

[22] S. P. Redmond, On zero divisor graph of small finite commutative rings, *Disc. Math.*, 307 (2007), 1155–1166.

[23] S. P. Redmond, An Ideal based zero divisor graph of a commutative ring, *Comm. Algebra*, 31 (9) (2003), 4425–4443.

[24] A. Sebo and E. Tannier, On metric generators of graphs, *Math. Oper. Research.*, 29 (2) (2004), 383–393.

[25] P. J. Slater, Dominating and reference sets in a graph, *J. Math, Phys. Sci.*, 22 (1988), 445–455.

[26] P.J. Slater, Leaves of trees, *Congressus Numerantium*, 14 (1975), 549– 559.

[27] S. Spiroff and C. Wickham, A zero divisor graph determined by equivalence classes of zero divisors, *Comm. Algebra*, 39 (7) (2011), 2338–2348.

[28] Simanjuntak, R., Siagian, P., & Vetrik, T. (2017). The multiset dimension of graphs. arXiv preprint arXiv:1711.00225.

[29] Pirzada, S., & Bhat, M. I. (2018). Computing metric dimension of compressed zero divisor graphs associated to rings. *Acta Universitatis Sapientiae, Mathematica*, *10*(2), 298-318.

[30] Kelenc, A., Tratnik, N., & Yero, I. G. (2018). Uniquely identifying the edges of a graph: the edge metric dimension. Discommutative ringete Applied Mathematics, 251, 204-220.

[31] Ali, N., Kousar, Z., Safdar, M., Tolasa, F. T., & Suleiman, E. (2023). Mapping Connectivity Patterns: Degree-Based Topological Indices of Corona Product Graphs. *Journal of Applied Mathematics*, *2023*.

[32] Safdar, M., Mushtaq, T., Ali, N., & Akgül, A. (2023). On study of flow features of hybrid nanofluid subjected to oscillatory disk. *International Journal of Modern Physics B*, 2450356.

[33] Mahboob, A., Hussain, T., Akram, M., Mahboob, S., Ali, N., & Raza, A. (2020). Characterizations of Chevalley groups using order of the finite groups. *Journal of Prime Research in Mathematics*, *16*(1), 46-51.